\documentclass[11pt]{amsart}

\usepackage[foot]{amsaddr}
\usepackage{amsmath}
\usepackage{amssymb}
\usepackage{amsthm}
\usepackage{bbm}
\usepackage[format=plain]{caption}
\usepackage{chngcntr}
\usepackage[margin=3cm]{geometry}
\usepackage{graphicx}
\usepackage{lmodern}
\usepackage{mathtools}
\usepackage{natbib}
\usepackage{subcaption}
\usepackage[table]{xcolor}

\makeatletter
\renewcommand\subsection{\@startsection{subsection}{2}
    \z@{-.5\linespacing\@plus-.7\linespacing}{.5\linespacing}
    {\centering\normalfont\scshape}
}
\def\@setfoot@addresses{
    {\itshape Institution}:
    \def\author##1{}%
    \def\\{\unskip, \ignorespaces}%
    \newif\if@firstaddr\@firstaddrtrue%
    \def\address##1##2{
    \if@firstaddr\@firstaddrfalse\else\par\fi\@ifnotempty{##1}{
        (\ignorespaces##1\unskip)}{\ttfamily\ignorespaces##2}
    }%
    \def\email##1##2{}%
    \def\curraddr##1##2{}%
    \def\urladdr##1##2{}%
    \addresses.
}
\makeatother

\newtheorem{ex}{Example}[section]
\newtheorem{lem}{Lemma}[section]
\newtheorem{pro}{Proposition}[section]
\newtheorem{rmk}{Remark}[section]
\newtheorem{thm}{Theorem}[section]
\numberwithin{figure}{section}

\begin{document} 

\address{Department of Mathematics, University of Michigan}
\author{Erhan Bayraktar and Thomas Bernhardt}
\date{\today}
\email{Erhan@umich.edu, Bernt@umich.edu}
\title{On the Continuity of the Root Barrier}

\maketitle

\begin{abstract}
\noindent We show that the barrier function in Root's solution to the Skorokhod embedding problem is continuous and finite at every point where the target measure has no atom and its absolutely continuous part is locally bounded away from zero.
\\[1em]
\end{abstract}

\section{Introduction}

\noindent In this paper, we consider the \cite{Skorokhod1961} embedding problem,  in particular the \cite{Root1969} solution in the formulation of \cite{Loynes1970}. Originally that meant, for given centered probability measure $\mu$ and given standard Brownian motion $W$, to find a lower semi-continuous function $\,r\colon[-\infty,\infty]\to[0,\infty]\,$ with $\,r(\pm\infty)=0\,$ such that $\,\tau=\inf\{t\geq0\,|\,t\geq r(W_t)\}\,$ fulfills $\,\mu\sim W_\tau\,$ and $\,(W_{t\wedge\tau})_{t\geq0}\,$ is uniformly integrable.

\cite{Hobson1998} observed that a continuous stock price, for which prices to vanilla options are known, can be connected to a Skorokhod embedding problem. More precisely, under the risk-neutral measure, the stock price is a time-changed Brownian motion and the prices of vanilla options determine the law of the stock price at maturity. In this setting, $\,(W_{t\wedge\tau})_{t\geq0}\,$ is a model for the time-changed price process and $\mu$ is the law at maturity. Note, here $\tau$ is a solution to the general Skorokhod embedding problem that may be different from Root's solution, i.e.\ $\tau$ may only be a stopping time such that $\,\mu\sim W_\tau\,$ and $\,(W_{t\wedge\tau})_{t\geq0}\,$ is uniformly integrable. \cite{Dupire2005} argued that, if $\tau$ itself represents a security of interest in this model, then the expected residual of Root's solution $\tau$ gives a lower bound on the call option of that security. He simply recalled that \cite{Rost1976} proved that Root's solution is of minimal expected residual, which can be stated using call option price, i.e.\ Root's solution minimizes $\,\mathbb{E}[(\tau-t)^+]\,$  for all $\,t>0\,$ under solutions $\tau$ to the Skorokhod embedding problem. \citet[Section 5]{CarrLee2010}, then, replaced the Brownian motion $W$ in the Skorokhod embedding problem with a geometric Brownian motion $S$ so that $\tau$ represents the terminal value of realized volatility. This established the importance of Root's solution to finding no-arbitrage lower bounds for VIX call options (VIX goes back to \cite{BrGa1989} and its call options are used to hedge short-term increases in volatility at the \cite{Cboe2019}). However, a rigorous proof that these were indeed the sharp bounds under all no arbitrage models were first given in \citet[Theorem 6.4]{CoWa2013}. 

Root gave an existence proof of his solution to the Skorokhod embedding problem but missed to give a construction. For a geometric Brownian motion, \citet[Section 5]{CarrLee2010} showed that the barrier function $r$ of Root's solution is the boundary of a free-boundary problem, and hence developed a numerical scheme to solve for $r$. First, solve the free-boundary problem and then identify the boundary, which is equal to $r$. That numerical scheme has been generalized to a wide class of processes $X$ beyond Brownian motion $W$ and geometric Brownian motion $S$. \cite{CoWa2013} generalized that approach to time-homogeneous martingale diffusion $X$ and \cite{GaObRe2015} to the time-inhomogeneous case. In their conclusion, Cox and Wang mentioned that it may be possible to say something about the regularity or shape of $r$ now, which has been challenging before, using the analytic literature on free boundary problems. Indeed, \cite{GaMiOb2015} followed that idea and proved that for Brownian motion $W$, the barrier function $r$ is continuous and hill shaped when $\mu$ is supported on a compact and its density is symmetric and V-shaped. But, to the best of our knowledge, there is nothing known about the regularity of $r$ beyond that. 

However, there is a result that relies on a priori knowledge about the continuity of $r$. In view of \cite{GaMiOb2015}, the barrier function $r$ solves a nonlinear Volterra integral equation of the first kind and conversely, that equation has a unique solution if the solution is continuous. In addition to continuity, if $r$ is hill shaped, then solving the Volterra equation is a fast way to compute numerically the barrier function $r$, in particular, computational simpler than solving a free boundary problem and finding the boundary function like described before. 

Unfortunately, in the current past, most progress related to Root's solution has been made in generalizing existence results instead of regularity results. For example, \citet[Section 7]{BeCoHu2017} extended the class of processes $X$ for which a Root solution exists to one dimensional regular Feller processes. They aimed to unify many known solutions to the Skorokhod embedding problem using ideas from optimal transport (for an overview of many solutions see \cite{Ob2004}). \cite{GaObZou2019} extended the free boundary problem approach to sufficiently regular Markov processes including discontinuous L\'{e}vy processes. On the other hand, \cite{CoObTo2019} and \cite{BeCoHu2020} generalized the formulation of the Skorokhod embedding problem with one target measure $\mu$ to a problem with a sequence of target measures $(\mu_i)_{i=1}^N$, and \cite{RiTaTo2020} even to a continuum of target measures $(\mu_i)_{i\in[0,1]}$. But, to the best of our knowledge, the literature lacks general regularity results on the barrier function $r$.

In this paper, we focus on the continuity of the barrier function $r$ of Root's solution. We make two observations that will explain what happens when $r$ is discontinuous and find a condition that will exclude those possibilities.
{
    \vskip 2em
    \begin{minipage}{0.49\textwidth}
        \centering
        \includegraphics[trim=0 1.5cm 0 3cm,clip,width=\linewidth]{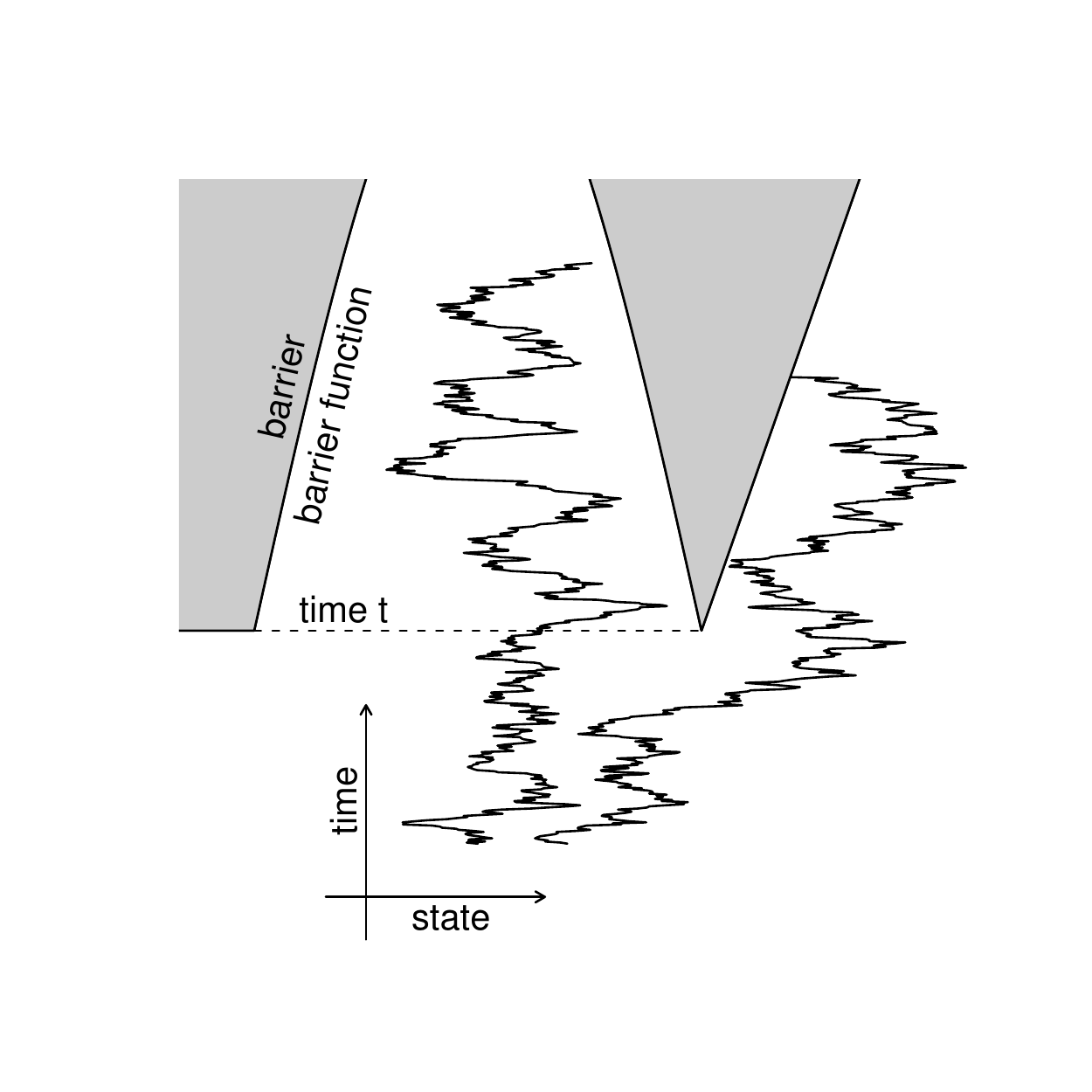}
    \end{minipage}
    \begin{minipage}{0.49\textwidth}
        \centering
        \includegraphics[trim=0 1.5cm 0 3cm,clip,width=\linewidth]{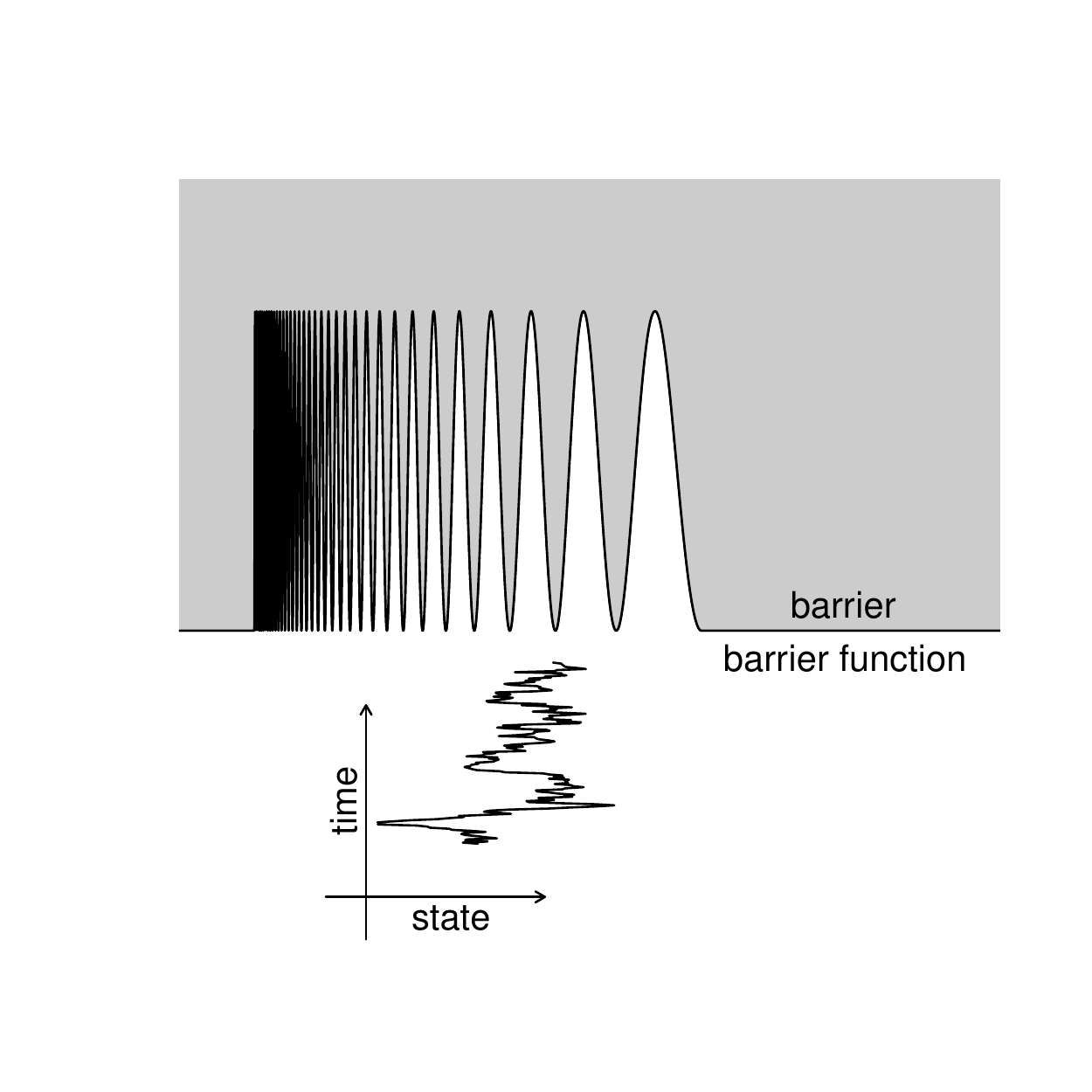}
    \end{minipage}
    \captionof{figure}{Two Brownian particles that are separated by a barrier after time $t$ (left). One Brownian particle that approaches a barrier with a discontinuous barrier function (right).}
    \label{fig:ideas}
    \vskip 1em
}

Figure \ref{fig:ideas}(left) shows two sample paths of Brownian motion with one sample path in a corridor with walls that are formed by the barrier function. The other sample path outside the corridor is unable to enter the corridor, as the path would be stopped at its walls. Hence, the mass of the target measure $\mu$ in the area of the corridor is at most the mass that passes through the entrance of the corridor. That phenomenon only depends on the continuity of the underlying process $W$ used in the Skorokhod problem.

Figure \ref{fig:ideas}(right) shows a Brownian sample path that approaches an essential discontinuity of the barrier function $r$. Note, jump discontinuities are directly linked to atoms in the target measure and cause no issue in our analysis. We can observe that an essential discontinuity in the barrier function leads to thin and thinner corridors close to the discontinuity. As the height of the corridors stays unchanged, we expect that a particle entering a thin corridor will be stopped early on at its walls and hardly reaches the end of the corridor. That argument relies on the volatility of the particle to stay away from zero over time and space. As the mass in the corridor is determined by the mass that passes through the entrance of the corridor, we expect little mass at the end of the corridor. Moreover, with thin and thinner corridors, we expect vanishing mass at the end of the corridors approaching the discontinuity. A complication in our analysis is that, in general, the mass that enters the thin and thinner corridors decreases with the width of the corridor to zero approaching the discontinuity. We will need to rescale the mass that arrives at the end of the corridors with the width of the corridors.

Overall, we expect a vanishing density of the target measure $\mu$ at points of the end of the corridors when approaching the discontinuity in the barrier function $r$. 

The paper is organized as follows before we start with our analysis, we specify the process and target measure that is used in the Skorokhod embedding problem and remark the precise assumptions that are indeed needed for our results in subsection \ref{subsec:Setting}. Then, we gather some notation that will be used throughout the paper in subsection \ref{subsec:Notation}. Section \ref{sec:Continuity} contains our analysis including Lemmata \ref{lem:continuous crossing}, \ref{lem:estimates} and \ref{lem:rzero} that formalize the heuristics of Figure \ref{fig:ideas}. Our main result can be found in Theorem \ref{thm:con and finite} and states that a target measure $\mu$ with a density that is locally bounded away from zero has a continuous barrier function $r$.

We finish the paper with an example resembling the situation in Figure \ref{fig:ideas}(right). 

\subsection{Setting}\label{subsec:Setting}

\noindent Denote by $(\Omega,\mathcal{F},\mathbb{P})$ a complete probability space that supports a standard Brownian motion $W$ and a random variable $X_0\sim\mu_0$ that are independent of each other. Let $\,\sigma\colon[0,\infty)\times\mathbb{R}\to\mathbb{R}\,$ such that there is $\,k\geq0\,$ with
\begin{gather}
    |\sigma(t,x)-\sigma(t,y)|\leq k|x-y|\quad\mbox{for all $\,t\geq0\,$ and $\,x,y\in\mathbb{R}$},\label{ass:1}
    \\|\sigma(t,x)|\leq k(1+|x|)\quad\mbox{for all $\,t\geq0\,$ and $\,x\in\mathbb{R}$.}\label{ass:2}
\end{gather}
The assumptions (\ref{ass:1}) and (\ref{ass:2}) ensure, see for example \citet[Theorem 17.2.3]{CoEll2015}, that there is a unique Markov process $X$ that is a strong solution to the dynamics
\[\mathrm{d}X_t=\sigma(t,X_t)\,\mathrm{d}W_t\quad\mbox{with initial condition $X_0$}.\]
Furthermore, assume that there is an open set $\mathcal{O}$ such that
\begin{gather}
    \mbox{$\sigma$ is bounded away from zero in $[0,\infty)\times K$ for any compact $\,K\subset\mathcal{O}$},\label{ass:3}
     \\\mbox{$\sigma$ is $\,\mathcal{C}^{1,2}$ in $[0,\infty)\times\mathcal{O}$},\label{ass:4}
    \\\mbox{$X,\mu_0$ take values in $\mathcal{O}$}.\label{ass:5}
\end{gather}
Lastly, let $\mu$ be a probability measure on $\mathcal{O}$ such that $\mu_0$ and $\mu$ are in convex order, i.e.\
\begin{gather}
    u_{\mu_0}(x)=-\int_\mathbb{R}|x-y|\,\mathrm{d}\mu_0(y)\geq -\int_\mathbb{R}|x-y|\,\mathrm{d}\mu(y)=u_\mu(x)>-\infty\quad\mbox{for all $\,x\in\mathbb{R}$}.\label{ass:6}
\end{gather}

The assumptions (\ref{ass:1}), (\ref{ass:2}), (\ref{ass:3}) and (\ref{ass:5}), (\ref{ass:6}) ensure that Root's solution to embed $\mu$ into $X$ exists, see \cite{GaObRe2015}, that means that there is a  lower semi-continuous function $\,r\colon[-\infty,\infty]\to[0,\infty]\,$ with $\,r(\pm\infty)=0\,$ such that $\,\tau=\inf\{t\geq0\,|\,t\geq r(X_t)\}\,$ fulfills $\,\mu\sim X_\tau\,$ and $\,(X_{t\wedge\tau})_{t\geq0}\,$ is uniformly integrable.

It is worth mentioning that the open set $\mathcal{O}$ has been introduced to the assumptions to specifically allow $X$ to be a geometric Brownian motion ($\sigma(t,x)=x\,$ and $\,\mathcal{O}=(0,\infty)$). Recall, Root's solution to the Skorokhod embedding problem with underlying geometric Brownian motion is linked to lower bounds on no-arbitrage prices of VIX call options, see \citet[Section 5]{CarrLee2010}.

Assumption (\ref{ass:4}) ensures that the marginals of the stochastic process $X$ have continuous densities, see \cite{Rogers1985}. That assumption is a technical requirement that allows us to bound the density of the absolutely continuous part of the target measure $\mu$ in terms of the marginal densities.  

Assumption (\ref{ass:6}) is necessary for $\,(X_{t\wedge\tau})_{t\geq0}\,$ to be uniformly integrable and the literature agrees on that this is the correct description to express that $\tau$ is ``small'', see \citet[Section 8]{Ob2004}. However, the assumption of uniform integrability in the Skorokhod embedding problem can be weakened, e.g.\ \cite{Monroe1972} and \cite{Wang2020}. 

Furthermore, we assume that the barrier function $r$ is regular, i.e.\
\begin{gather}\label{ass:7}
    u_{\mu_0}(x)=u_\mu(x)\quad\mbox{implies}\quad r(x)=0\quad\mbox{for all $\,x\in\mathbb{R}$}.
\end{gather}
\citet{GaObRe2015} introduced that definition of regularity that ensures uniqueness of the existing barrier function. It is a generalization of \citet{Loynes1970} regularity that was tailored to the special case $\,\mu_0=\delta_0$, the Dirac delta measure at zero, with the same purpose of uniqueness. Loynes defined the interval $\,(\alpha,\beta)\,$ with $\,\alpha=\sup\{x<0\,|\,r(x)=0\}\,$ and $\,\beta=\inf\{x>0\,|\,r(x)=0\}$, and defined that $r$ is regular if it vanishes outside of $\,(\alpha,\beta)$. In the same way, we can restate Gassiat et al.'s regularity criterion as $r$ is regular if it vanishes outside of 
\[\mathcal{I}=\{x\,|\,u_{\mu_0}(x)\neq u_\mu(x)\}.\]

\begin{rmk}
    {\normalfont Even though, we assume {\normalfont(\ref{ass:1})-(\ref{ass:7})} to fit in the current literature on Root's solution, our analysis only requires the assumptions: {\normalfont(i)} there is a Root solution $r$, {\normalfont(ii)} $r$ is minimal in the sense that, if $q$ is a barrier function solving the Skorokhod embedding problem with $\,q\leq r$, then $\,q=r$, {\normalfont(iii)} the underlying process $X$ is a continuous Markov process with continuous marginal densities, and {\normalfont(iv)} all two-sided hitting times of $X$ in $\mathcal{O}$ hit both barriers after any time with positive probability, i.e.\ for any $\,(a,b)\subset\mathcal{O}\,$ and $\,t_2>t_1\geq s>0\,$ the hitting time $\,\varsigma=\inf\{t\geq0\,|\,X_t\notin(a,b)\}\,$ fulfills
    \begin{gather*}
        \mathbb{P}[X_\varsigma=a,\varsigma\in(t_1,t_2)\,|\,X_s\in(a,b)]>0,
        \\\mathbb{P}[X_\varsigma=b,\varsigma\in(t_1,t_2)\,|\,X_s\in(a,b)]>0.
    \end{gather*}
    }
\end{rmk}

\subsection{Notation}\label{subsec:Notation}

\noindent To distinguish easily between functions and measures, we will consistently denote evaluations of functions using parentheses, e.g.\ $\,r(x)\,$ for some $\,x\in\mathbb{R}$, and evaluations of measures using brackets, e.g.\ $\,\mu[(x,y)]\,$ for some $\,x<y$. We will allow functions to operate on sets, i.e.\ the elementwise evaluation 
\begin{gather*}
    r(A)=\{r(a)\,|\,a\in A\}\quad\mbox{for all $\,A\subset\mathbb{R}$},
    \\r(A)\geq t\quad\mbox{if and only if}\quad r(a)\geq t\quad\mbox{for all $\,a\in A$},\quad\mbox{and given $\,t\geq0$}. 
\end{gather*}
Most notably, $\,r(A)=\infty\,$ means $\,r(a)=\infty\,$ for all $\,a\in A$. Moreover, we will express our results using the following measures,
\begin{align*}
    \mu_t[A]=\mathbb{P}[X_{t\wedge\tau}\in A]\quad\mbox{for $\,t\geq0\,$ and Borel set $A$}.
\end{align*}

\section{Continuity}\label{sec:Continuity}
The following lemmata formalize the heuristics of Figure \ref{fig:ideas}. Lemma \ref{lem:continuous crossing} captures the phenomoenon of Figure \ref{fig:ideas}(left). Lemma \ref{lem:estimates} contains an estimate of the mass that arrives at the end of a corridor described in Figure \ref{fig:ideas}(right). And Lemma \ref{lem:rzero} states that a discontinuity in the barrier function $r$ leads to a sequence of points where an approximate density of $\mu$ decreases to zero.  
\begin{lem}\label{lem:continuous crossing}If $\,t\geq s\geq r(x)\vee r(y)\,$ for some  $\,y>x$, then
    \begin{itemize}
        \item[(i)]$\{X_{\varsigma}\in(x,y),\,\tau\geq\varsigma\geq t\}\subset\{X_s\in(x,y)\}\,$ for any (possibly random) time $\varsigma$,
        \item[(ii)]$\,\mu[(x,y)]\leq\mu_t[(x,y)]$.
    \end{itemize}
\begin{proof}
    If $\,X_\varsigma\in (x,y)$ and $\,X_s\notin(x,y)$ and $\,\varsigma\geq t\geq s$, then continuity of $X$ ensures that there is a time $\,\varsigma>\hat{s}\geq s\,$ with $\,X_{\hat{s}}\in\{x,y\}$. This implies, as $\,\hat{s}\geq s\geq r(x)\vee r(y)$, that $\,\hat{s}\geq r(X_{\hat{s}})$. In view of the definition of $\tau$,
    this means that $\,\hat{s}\geq\tau\,$ and therefore $\,\varsigma>\tau$. Hence, (i).  
     
    Now, (ii) follows from (i), as $\,\varsigma=\tau\,$ and $\,s=t\,$ reveals $\,\{X_\tau\in(x,y),\tau\geq t\}\subset\{X_t\in(x,y)\}$, in particular,
    \begin{align*}
        \mu[(x,y)]&=\mathbb{P}[X_\tau\in (x,y)]
        \\&=\mathbb{P}[X_\tau\in (x,y),\,t>\tau]+\mathbb{P}[X_\tau\in (x,y),\,t\leq\tau]
        \\&\leq\mathbb{P}[X_\tau\in (x,y),\,t>\tau]+\mathbb{P}[X_t\in (x,y),\,t\leq\tau]
         \\&=\mathbb{P}[X_{t\wedge\tau}\in(x,y)]
        \\&=\mu_t[(x,y)].
    \end{align*}
\end{proof}
\end{lem}

\begin{lem}\label{lem:estimates} 
    If $\,r(A)\geq t\geq s\geq r(x)\vee r(y)\,$ for a Borel subset $A$ of $(x,y)$ with $\,[x,y]\subset\mathcal{O}$, and $\lambda$ denotes the Lebesgue measure, then there is a constant $k_{x,y}$ such that
    \[\mu_t[A]\leq\lambda[A]\,k_{x,y}\,\mu_s[(x,y)].\]
    Moreover, the constant can be chosen so that $\,k_{\hat{x},\hat{y}}\leq k_{x,y}\,$ whenever $\,(\hat{x},\hat{y})\subset(x,y)$.
\begin{proof}
    Before we begin with the estimation, we make the following observations:
    
    1) $\,\{X_\tau\in A,\,t>\tau\}\,$ is an empty set. To see this, observe that $\,X_\tau\in A\,$ implies $\,r(X_\tau)\in r(A)$. But on the other hand, if $\,t>\tau$, then $\,r(A)\geq t>\tau\geq r(X_\tau)$, and therefore $\,r(X_\tau)\notin r(A)$. Note that $\,r(A)\geq t\,$ is an assumption and $\,\tau\geq r(X_\tau)$ is a consequence of $r$ being lower semi-continuous.
    
    2) $\,\{X_t\in A,\,\tau\geq t\}=\{X_t\in A,\,\tau\geq t,\,X_s\in(x,y)\}$, which is a consequence of Lemma \ref{lem:continuous crossing}(i) following from setting $\,\varsigma=t$.
    
    3) Let $\mathcal{F}_s$ be the $\sigma$-algebra generated by the process $X$ up to to time $s$, then the Markov property of $X$, see for example \citet[Theorem 17.2.3]{CoEll2015}, yields $\,\mathbb{P}[X_t\in A\,|\,\mathcal{F}_s]=\mathbb{P}[X_t\in A\,|\,X_s]$ almost surely.
    
    4) In view of \cite{Rogers1985}, there is a continuous function $\,f\colon\mathbb{R}\times\mathbb{R}\to[0,\infty)\,$ so that for all Borel sets $\,B\subset[x,y]\,$ the equation $\,\mathbb{P}[X_t\in B\,|\,X_s]\mathbbm{1}_{X_s\in[x,y]}=\int_B f(b,X_s)\,\mathrm{d}b\,\mathbbm{1}_{X_s\in[x,y]}\,$ holds true $\mathbb{P}$-almost surely. As $\,[x,y]\times[x,y]$ is compact, there is a constant $\,\infty>k_{x,y}\geq f\,$ on that compact. 
    
    \noindent Having these four observations in mind,
    \begin{align*}
        \mu_t[A]&=\mathbb{P}[X_{t\wedge\tau}\in A]
        \\&=\mathbb{P}[X_t\in A,\,t\leq\tau]+\mathbb{P}[X_\tau\in A,\,t>\tau]
        \\&\stackrel{1)}{=}\mathbb{P}[X_t\in A,\,t\leq\tau]
        \\&\stackrel{2)}{=}\mathbb{P}[X_t\in A,\,t\leq\tau,\,X_s\in(x,y)]
        \\&\leq\mathbb{P}[X_t\in A,\,s\leq\tau,\,X_s\in(x,y)]
        \\&=\mathbb{E}\big[\mathbb{P}\big[X_t\in A\,\big|\,\mathcal{F}_s\big]\mathbbm{1}_{X_s\in(x,y)}\mathbbm{1}_{s\leq\tau}\big]
        \\&\stackrel{3)}{=}\mathbb{E}\big[\mathbb{P}\big[X_t\in A\,\big|\,X_s\big]\mathbbm{1}_{X_s\in(x,y)}\mathbbm{1}_{s\leq\tau}\big]
        \\&\stackrel{4)}{\leq}\mathbb{E}\Big[\lambda[A]\,k_{x,y}\,\mathbbm{1}_{X_s\in(x,y)}\mathbbm{1}_{s\leq\tau}\Big]
        \\&\leq \lambda[A]\,k_{x,y}\,\mathbb{P}[X_{s\wedge\tau}\in(x,y)]
        \\&=\lambda[A]\,k_{x,y}\,\mu_s[(x,y)]).
    \end{align*}
    
    For the second part of the statement, it is enough to observe that $\,[\hat{x},\hat{y}]\subset[x,y]$, hence $\,k_{x,y}\geq f\,$ also on $\,[\hat{x},\hat{y}]\times[\hat{x},\hat{y}]$.
\end{proof}
\end{lem}

\begin{lem}\label{lem:rzero}
    If $\,\limsup_{y\downarrow x}r(y)>\liminf_{y\downarrow x}r(y)\,$ for some $x\in \mathcal{O}$, then 
    \[\liminf_{\substack{(\varepsilon,y)\rightarrow(0,x)\\\varepsilon>0,\;y>x}}\frac{\mu[(y,y+\varepsilon)]}{\varepsilon}=0.\]
\begin{proof}
    Choose times $\,t>s\,$ such that $\,\limsup_{y\downarrow x}r(y)>t\,$ and $\,s>\liminf_{y\downarrow x}r(y)$. Let $\,x_n\downarrow x\,$ such that $\,t< r(x_n)\rightarrow\limsup_{y\downarrow x}r(y)$. Consider
    \begin{align*}
        \delta_1^n&=\inf\{\delta>0\,|\,r(x_n-\delta)\leq t\},
        \\\delta_2^n&=\inf\{\delta>0\,|\,r(x_n+\delta)\leq t\}.
    \end{align*}
    By lower semi-continuity of $r$, it holds that $\,r(x_n-\delta_1^n)\leq t\,$ and  $\,r(x_n+\delta_2^n)\leq t$. As $\,r(x_n)>t$, it holds true that $\,\delta_1^n,\delta_2^n>0$. Let
    \[A_n=(x_n-\delta_1^n,x_n+\delta_2^n).\]
    Furthermore, choose $\,\hat{x}_n\downarrow x\,$ such that $\,s> r(\hat{x}_n)\rightarrow\liminf_{y\downarrow x}r(y)\,$ and $\,\hat{x}_n>x_n\,$ for all large enough $n$. Then, by construction, the set $\,A_n\subset(x,\hat{x}_n)\,$ fulfills the assumptions of Lemma \ref{lem:estimates}. 
    
    Note, $\,A_n\subset(x,\hat{x}_n)\,$ and $\,\hat{x}_n\downarrow x$ imply $\,\delta_1^n+\delta_2^n\rightarrow0$, in particular $\,\delta_1^n+\delta_2^n<\infty$. Now, putting Lemma \ref{lem:continuous crossing}(ii) and \ref{lem:estimates} together and using continuity of probability measures,
    \begin{align*}
        \frac{\mu[A_n]}{\delta_1^n+\delta_2^n}
        \leq\frac{\mu_t[A_n]}{\delta_1^n+\delta_2^n}
        \leq k_{x,\hat{x}_0}\,\mu_s[(x,\hat{x}_n)]\xrightarrow[]{n\uparrow\infty}0.
    \end{align*}
    As $\,\delta_1^n+\delta_2^n\rightarrow0$, the statement follows from considering the left endpoint $y_n$ and the length $\varepsilon_n$ of $\,A_n=(x_n-\delta_1^n,x_n+\delta_2^n)$.
\end{proof}  
\end{lem}

\begin{rmk}\label{rem:rzero}
    {\normalfont Lemma \ref{lem:rzero} can be slightly strengthened by keeping the same assumptions but requiring that the sequences in the conclusion not only fulfill $\,\varepsilon\downarrow0,\,y\downarrow x\,$ but also fulfill $\,(y,y+\varepsilon)\subset\mathcal{I}=\{x\,|\,u_{\mu_0}(x)\neq u_\mu(x)\}$. To see this, observe that $\,A_n=(x_n-\delta_1^n,x_n+\delta_2^n)\,$ from the proof fulfills $\,r(A_n)>t>s\geq0$, in particular, $A_n$ cannot be outside of $\,\mathcal{I}$ where $r$ is vanishing by regularity, hence $\,A_n\subset\mathcal{I}$.
    }
\end{rmk}

As the barrier function $r$ is lower semi-continuous, $r$ will always approach infinity continuously. Given a point $x$ where $r$ is infinite, but all other points are finite, then the barrier function describes a thinner and thinner corridor for the approaching process over time, which is the main argument in the previous lemma. So it comes to no surprise that the previous ideas can be used to analyze points where the barrier function is infinity.

\begin{lem}\label{lem:rinfzero}
    If $\,r(x)=\infty$, then
    \[\liminf_{\substack{(\varepsilon,y)\rightarrow (0,x)\\\varepsilon>0}}\frac{\mu[(y-\varepsilon,y+\varepsilon)]}{\varepsilon}=0.\]
\begin{proof}
   Here, we employ an alteration of the arguments in the proof of Lemma \ref{lem:rzero}. But first observe that $\,x\in\mathcal{O}$, because otherwise $X$ never attains $x$, hence $X$ never passes through $x$ by continuity, hence $r(x)$ can be chosen to be zero, in contradiction to $\,r(x)=\infty\,$ and the uniqueness of $r$.
   
   Look at
    \begin{align*}
        \delta_1^n&=\inf\{\delta>0\,|\,r(x-\delta)\leq n\},
        \\\delta_2^n&=\inf\{\delta>0\,|\,r(x+\delta)\leq n\}.
    \end{align*}
    
    First, consider the case $\,\delta_1^n+\delta_2^n\rightarrow0\,$ and observe:
    
    1) As $\,\delta_1^n+\delta_2^n\rightarrow0$, there is an integer $m$ such that $\,\delta_1^m+\delta_2^m<\infty$. Let $\,k=k_{x-\delta_1^m,x+\delta_2^m}$ from Lemma \ref{lem:estimates}.
    
    2) $\,A_n=(x-\delta_1^n,x+\delta_2^n)\subset(x-\delta_1^{n-1},x+\delta_2^{n-1})=A_{n-1}\,$ fulfills the assumptions of Lemma \ref{lem:continuous crossing}(ii) and \ref{lem:estimates}. 
    
    3) $\,A_n\downarrow A_\infty\,$ with $\,r(A_\infty)=\infty$. That convergence together with $\,X_{n-1\wedge\tau}\rightarrow X_\tau\,$ imply $\,\limsup_{n\rightarrow\infty}\mu_{n-1}[A_{n-1}]\leq\mu[A_\infty]$.
    
    4) As $r$ is lower semi-continuous, $\,\tau\geq r(X_\tau)$, hence $\,\{X_\tau\in A_\infty\}\subset\{\tau=\infty\}$, which is a $\mathbb{P}$-zero set.
    
    \noindent Altogether,
    \begin{align*}
        \limsup_{n\rightarrow\infty} \frac{\mu[A_n]}{\delta_1^n+\delta_2^n}
        &\stackrel{1),2)}{\leq}\limsup_{n\rightarrow\infty}\,k\,\mu_{n-1}[A_{n-1}]
        \\&\stackrel{3)}{\leq} k\,\mu[A_\infty]
        \\&=k\,\mathbb{P}[X_\tau\in A_\infty]
        \\&\stackrel{4)}{\leq} k\,\mathbb{P}[\tau=\infty]
        \\&=0.
    \end{align*}
    Keep in mind that $\,\delta_1^n+\delta_2^n\rightarrow0$, hence the statement follows from considering the middle point $y_n$ and its distance $\varepsilon_n$ to the boundary of $\,A_n=(x-\delta_1^n,x+\delta_2^n)$.
     
    Second, consider $\,\delta_1^n+\delta_2^n\not\rightarrow0$. Then without loss of generality, there is $\,\delta>0\,$ with $\,r((x-\delta,x))=\infty$. Hence, $\,\mu[(x-\delta,x)]=0$, and the statement follows from choosing a sequence $\,y_n\uparrow x\,$ such that $\,\delta/2>\varepsilon_n\colon\!\!\!\!=x-y_n>0\,$ and observing that $\,\mu[(y_n-\varepsilon_n,y_n+\varepsilon_n)]/\varepsilon_n=0/\varepsilon_n=0$.
\end{proof}  
\end{lem}

It is important to notice that Lemma \ref{lem:rzero} makes no statement about what happens at the point $x$ with the barrier function $r$. The next lemma shows that $r(x)$ can only differ from its limit inferior in the presence of an atom.

\begin{lem}\label{lem:ratom}
    If $\,\liminf_{y\downarrow x} r(y)>r(x)\,$ for some $\,x\in\mathcal{O}$, then $\mu$ has an atom at $x$.
\begin{proof} 
    Let 
    \[z=\sup\{u>x\,|\,r(y)>r(x)\;\;\mbox{for all $\,y\in(x,u)$}\}.\] 
    Using proof by contradiction, we show that $\,\mu_{r(x)}[(x,z)]>0$, i.e.\ assume $\,\mu_{r(x)}[(x,z)]=0$. 
    
    Define two new functions $\bar{r},\hat{r}$ with $\,\bar{r}(y)=\hat{r}(y)=r(y)\,$ for $\,y\notin(x,z)$ and $\,\bar{r}(y)=r(x)\,$ and $\,\hat{r}(y)=K\,$ for all $\,y\in(x,z)$ and any constant $\,K>r(x)$. Consider the corresponding $\,\bar{\tau}=\inf\{t\geq0\,|\,t\geq\bar{r}(X_t)\}$, $\,\hat{\tau}=\inf\{t\geq0\,|\,t\geq\hat{r}(X_t)\}\,$ and $\,\bar{\mu}=\mathbb{P}_{X_{\bar{\tau}}}$, $\hat{\mu}=\mathbb{P}_{X_{\hat{\tau}}}$, the image measures of $X_{\bar{\tau}}$, $X_{\hat{\tau}}$. Observe:
    
    1) As $\,r\geq\bar{r}$, it follows that $\,\tau\geq\bar{\tau}$, hence Lemma \ref{lem:continuous crossing}(i) implies $\,\{X_{\bar{\tau}}\in(x,z),\,\bar{\tau}\geq r(x)\}\subset\{X_{r(x)}\in(x,z)\}$.
    
    2) Consider an event such that $\,X_{\hat{\tau}}\neq X_{\bar{\tau}}$. As $\,r(x)\wedge\hat{\tau}=r(x)\wedge\bar{\tau}$, it follows that $\,\bar{\tau}\geq r(x)\,$ for that event. Moreover, $\,\hat{r}\geq\bar{r}\,$ implies $\,\hat{\tau}\geq\bar{\tau}\,$ like before, in particular, $\,\hat{\tau}>\bar{\tau}$. By definition of $\hat{\tau}$ and $\bar{\tau}$, note that $\bar{r}$ is lower semi-continuous, $\,\hat{r}(X_{\bar{\tau}})>\bar{\tau}\geq\bar{r}(X_{\bar{\tau}})$. By construction $\hat{r}$ and $\bar{r}$ only differ on $(x,z)$, hence $\,X_{\bar{\tau}}\in(x,z)$. Overall, $\,\{X_{\hat{\tau}}\neq X_{\bar{\tau}}\}\subset\{X_{\bar{\tau}}\in(x,z),\,\bar{\tau}\geq r(x)\}$.
    
    \noindent Having these two observations in mind and the assumption $\,\mu_{r(x)}[(x,z)]=0$,
    \begin{align*}
        0&=\mu_{r(x)}[(x,z)]
        \\&\stackrel{1)}{\geq}\mathbb{P}[X_{\bar{\tau}}\in(x,z),\,\bar{\tau}\geq r(x)]
        \\&\stackrel{2)}{\geq}\mathbb{P}[X_{\hat{\tau}}\neq X_{\bar{\tau}}].
    \end{align*}
    Note, in the above derivation we could have used $\,t\wedge\tau,t\wedge\bar{\tau}$ or $\,t\wedge\hat{\tau},t\wedge\bar{\tau}$ instead of $\hat{\tau},\bar{\tau}$, hence $\,\mathbb{P}[X_{t\wedge\tau}\neq X_{t\wedge\bar{\tau}}]=0\,$ and $\,\mathbb{P}[X_{t\wedge\hat{\tau}}\neq X_{t\wedge\bar{\tau}}]=0\,$ as well. Therefore, the uniform integrability of $\,t\mapsto X_{t\wedge\tau}$ implies that $\,t\mapsto X_{t\wedge\bar{\tau}}$ and $\,t\mapsto X_{t\wedge\hat{\tau}}$ are uniformly integrable. Also, in the above derivation we could have used $\,\tau,\bar{\tau}$ instead of $\hat{\tau},\bar{\tau}$, hence $\,\mathbb{P}[X_\tau\neq X_{\bar{\tau}}]=0\,$ and $\,\mathbb{P}[X_{\hat{\tau}}\neq X_{\bar{\tau}}]=0$, in particular, $\,\hat{\mu}=\mu=\bar{\mu}$. Overall, $\hat{r}$ and $\bar{r}$ are different barrier functions of Root solutions to embed $\mu$ into $X$ with initial measure $\mu_0$. Because Root solutions are unique on $\mathcal{I}$, but $\hat{r}$ and $\bar{r}$ differ on $(x,z)$, the whole interval $(x,z)$ does not belong to $\mathcal{I}$. Thus, by regularity of $r$, it holds that $\,r(y)=0\,$ for $\,y\in(x,z)$. By lower semi-continuity of $r$, it follows $\,r(y)=0\,$ for $\,y\in[x,z)\,$ in contradiction to $\,\liminf_{y\downarrow x} r(y)>r(x)$. So, our original assumption is wrong and $\,\mu_{r(x)}[(x,z)]>0\,$ is correct.
    
    Choose $\,x<\hat{z}<z\,$ with $\,\mu_{r(x)}[(x,\hat{z})]>0$. Due to maximality of $z$ together with lower semi-continuity of $r$ and $\,\liminf_{y\downarrow x}r(y)>r(x)$, it holds true that $\,\inf_{y\in(x,\hat{z})}r(y)>r(x)$. Hence, $\,\mu_{r(x)}[(x,\hat{z})]>0\,$ and $\,(x,\hat{z})\times[r(x),\inf_{y\in(x,\hat{z})}r(y))\neq\emptyset\,$ is free from the graph of the barrier function. As all two-sided hitting times of $X$ in $\mathcal{O}$ hit both barriers with positive probability in any time, that yields that $x$ is hit with positive probability, hence $x$ is an atom of $\mu$.    
\end{proof}
\end{lem}

Now, we can state our main continuity result. The previous lemmata show how a discontinuity in the barrier function leads to a vanishing estimate on the derivative of the target measure $\mu$. So, the natural thing to do is to exclude a vanishing density and in turn, make discontinuities impossible.  

\begin{pro}\label{pro:con}
    Consider a point $x$ that is not an atom of $\mu$.
    \begin{itemize}
        \item[(i)] If the absolutely continuous part of $\mu$ has a derivative that is bounded away from zero in $\,[x,y)\,\cap\,\mathcal{I}$, then $r$ is right-continuous at $x$,
        \item[(ii)] if the absolutely continuous part of $\mu$ has a derivative that is bounded away from zero in an open area around $x$, then $r$ is finite at $x$.
    \end{itemize} 
\begin{proof}
    Let $\mu^c$ be the absolutely continuous part of the measure $\mu$. Denote by $\mu'$ the density of $\mu^c$ and let $\,\mu'(z)\geq k>0\,$ for all $\,z\in[x,y)\,\cap\,\mathcal{I}\,$ by assumption. Then any sequence $\,\varepsilon_n\downarrow0\,$ and $\,y_n\downarrow x\,$ with $\,(y_n,y_n+\varepsilon_n)\subset\mathcal{I}\,$
    fulfills that
    \[\frac{\mu[(y_n,y_n+\varepsilon_n)]}{\varepsilon_n}\geq\frac{\mu^c[(y_n,y_n+\varepsilon_n)]}{\varepsilon_n}=\frac{\int_{y_n}^{y_n+\varepsilon_n}\!\mu'(z)\,\mathrm{d}z}{\varepsilon_n}\geq\frac{k\,\varepsilon_n}{\varepsilon_n}=k>0.\]
    Hence, Remark \ref{rem:rzero} implies that $\,\limsup_{y\downarrow x}r(y)=\liminf_{y\downarrow x}r(y)$. As there is no atom at $x$, Lemma \ref{lem:ratom} implies $\,r(x)=\liminf_{y\downarrow x}r(y)$, hence, $r$ is right-continuous at $x$ and (i) has been proven.
    
    In a similar way, a density bounded away from zero around $x$ ensures that the consequence of Lemma \ref{lem:rinfzero} cannot hold, i.e.\ $r(x)<\infty$.
\end{proof}  
\end{pro}

Now, we are ready to formulate our main result.

\begin{thm}\label{thm:con and finite}
    If $\mu$ is atom-free and its absolutely continuous part has a derivative that is locally bounded away from zero when restricted to $\,\mathcal{I}=\{x\,|\,u_{\mu_0}(x)\neq u_\mu(x)\}$, then $r$ is continuous and finite.
\begin{proof}
    As $\,x\mapsto u_\nu(x)=-\int_\mathbb{R}|x-y|\,\mathrm{d}\nu(y)\,$ for $\,\nu\in\{\mu_0,\mu_1\}\,$ is continuous, $\,\mathcal{I}=\{x\,|\,u_{\mu_0}(x)\neq u_\mu(x)\}\,$ is an open set. Hence, for each $\,x\in\mathcal{I}$, there is 
    an open set that contains $x$ and in which the absolutely continuous part of $\mu$ is bounded away from zero. Therefore, Proposition \ref{pro:con}(ii) ensures that $r$ is finite in $\,\mathcal{I}$. As $r$ vanishes outside of $\,\mathcal{I}$ by regularity, $r$ is finite everywhere. 

    As the density of $\mu$ is locally bounded away from zero in $\,\mathcal{I}\,$ by assumption, Proposition \ref{pro:con}(i) ensures that $r$ is right-continuous. Furthermore, looking at $\,\hat{\nu}_0[A]\colon\!\!\!\!=\nu[-A]\,$ for Borel set $A$ and $\,\nu\in\{\mu_0,\mu\}\,$ as well as $\,\hat{r}(x)=r(-x)\,$ for $\,x\in[-\infty,\infty]\,$ and applying Proposition \ref{pro:con}(i) to it as before shows that $r$ is left-continuous.
\end{proof}  
\end{thm}

\begin{rmk}\label{rmk:application}
    \textnormal{
    \cite{AnStr2011}, \cite{AnHoStr2015} and \cite{AnEnFrReis2020} analyzed Bass' solution to the Skorokhod embedding problem. One simple consequence of their approach is that an embeddable measure $\mu$ with a density that is bounded away from zero in its compact support, can be embedded in bounded time. Considering Root's solution and Theorem \ref{thm:con and finite}, that simple result is the consequence of a continuous function achieving its maximum on a compact support.
    }
    
    \textnormal{
    Moreover, \cite{AnStr2011} and \cite{AnHoStr2015} obtained a condition on the local behaviour of $\mu$ that is necessary so that $\mu$ can be embedded in bounded time. Assuming $\mu$ has a density, that condition makes only restrictions on the points where the density vanishes and makes no assertions for the points when the density is strictly positive. One may ask if any condition on the local behaviour of $\mu$ at points with strictly positive density is needed that is just not captured by their condition. Considering Root's solution and Theorem \ref{thm:con and finite}, the answer is no; if the density is locally away from zero then the barrier function is finite and continuous and hence locally bounded.
    }
\end{rmk}

We conclude the paper with an example that resembles the situation in Figure \ref{fig:ideas}(right) and show that conditions on the density are, indeed, needed to ensure continuity of the barrier function $r$. To our knowledge, a direct way to construct a target measure from a barrier function or vice versa with arbitrary properties is unknown. Our approach is implicit and relies on general existence results when $X$ is a local martingale. We start with an observation about regular barrier functions that will help us to paste different barrier functions of different measures together and the result will be the barrier function of the sum of these measures.

\begin{lem}\label{lem:mu[tail]=0}
    If $X$ is a local martingale and $\,\mu[(x,\infty)]=0$, then $\,r(y)=0\,$ for all $\,y\geq x$.
\begin{proof}
    Using proof by contradiction, we show that $\,\mu_0[(x,\infty)]=0$, so we assume $\,\mu_0[(x,\infty)]>0$.   
    
    As $\,t\mapsto X_{t\wedge\tau}\,$ is a uniformly integrable martingale, $\,\mathbb{E}[X_0]=\mathbb{E}[X_\tau]$, in particular, as $\,\mu[(x,\infty)]=0$, it holds that $\,\mathbb{E}[X_\tau]\leq x\,$ and in turn $\,\mathbb{E}[X_0]\leq x\,$ and therefore $\,\mu_0[(-\infty,x)]>0$, too. Hence, Jensen's inequality applied to $\mu_0$ and the convex function $\,y\mapsto|x-y|$ is strict, i.e.
    \begin{align*}
        \int_\mathbb{R}\!|x-y|\,\mathrm{d}\mu_0(y)
        &>\Big|x-\!\int_\mathbb{R}\!y\,\mathrm{d}\mu_0(y)\Big|
        \\&=\big|x-\mathbb{E}[X_0]\big|
        \\&=\big|x-\mathbb{E}[X_\tau]\big|
        \\&=\Big|x-\!\int_{(-\infty,x]}\!y\,\mathrm{d}\mu(y)\Big|
        \\&=\int_\mathbb{R}\!|x-y|\,\mathrm{d}\mu(y).
    \end{align*}
    That is a contradiction to the assumed convex order between $\mu_0$ and $\mu$. So our original assumption is wrong and $\,\mu_0[(x,\infty)]=0\,$ is correct.
    
    Now, the statement follows from regularity of $r$, as $\,\mu_0[(x,\infty)]=\mu[(x,\infty)]=0\,$ and $\,\mathbb{E}[X_0]=\mathbb{E}[X_\tau]\,$ imply $\,u_{\mu_0}(y)=u_{\mu}(y)\,$ for all $\,y\geq x$.
\end{proof}
\end{lem}

As mentioned before, we conclude the paper with an example that resembles the situation in Figure \ref{fig:ideas}(right). We construct probability measures with barrier functions that are zero outside some interval and are infinity at some open neighborhood inside the interval. Pasting those barriers with thinner and thinner intervals together yields two sequences, for which the first sequence evaluated at the barrier function is bounded and the second sequence is infinity. In particular, the barrier function cannot be continuous. More precisely:

\begin{ex}
    Let $\,X=W$ a standard Brownian motion with $\,\mu_0=\delta_0$, the Dirac delta measure. Then, there is $\mu$ with a density that is not bounded away from zero and $r$ is discontinuous. 
\begin{proof}
    For given interval $(a,b)$, consider the probability measure $\xi_{\lambda,p}$ for $\,\lambda\in[0,1]$ and $\,p>0\,$ with density
    \[\xi_{\lambda,p}'(x)=(1-\lambda)\,\frac{2(a+p-x)}{p^2}\mathbbm{1}_{x<a+p}+\lambda \,\frac{2(x-b+p)}{p^2}\mathbbm{1}_{x>b-p}\quad\mbox{for $\,x\in(a,b)$}.\]
    For any $c\in(a,b)$, there are $\,\lambda^*\in(0,1)\,$ and $\,p^*\in(0,\frac{b-a}{4})\,$ such that the first moment of $\xi_{\lambda^*,p^*}$ is $c$. To see this, observe that $\,\xi_{0,p}\Rightarrow\delta_{a}$, the Dirac measure at $a$, and $\,\xi_{1,p}\Rightarrow\delta_{b}$, the Dirac measure at $b$, when $\,p\downarrow0$, and observe that the first moment of $\xi_{\lambda,p}$ is continuous with respect to $\lambda$ and $p$. Now, choose $p^*$ so that $c$ is strictly in-between the first moments of $\xi_{0,p^*}$ and $\xi_{1,p^*}$, and then, choose $\lambda^*$ so that the first moment of $\xi_{\lambda^*,p^*}$ equals $c$. For ease of notation, we write
    \[\xi_c=\xi_{\lambda^*,p^*}.\]
    Define, for a given finite measure $\nu$ on $(a,b)$,
    \[\eta[A]=\int_a^b\xi_c[A]\,\mathrm{d}\nu\quad\mbox{for Borel set $A$}.\]
    Observe:
    
    1) $\eta$ is absolutely continuous, because any Lebesgue zero-set is a zero-set for all $\xi_c$.
    
    2) $\eta[\mathbb{R}]=\nu[\mathbb{R}]$. Moreover, Jensen's inequality yields for $\,x\in\mathbb{R}\,$ that
    \begin{align*}
        \int_a^b|x-y|\,\mathrm{d}\eta(y)
        &=\int_a^b\int_a^b|x-y|\,\mathrm{d}\xi_c(y)\,\mathrm{d}\nu(c)
        \\&\geq\int_a^b\Big|\,x-\!\int_a^b\!y\,\mathrm{d}\xi_c(y)\Big|\,\mathrm{d}\nu(c)
        \\&=\int_a^b|x-c|\,\mathrm{d}\nu(c).
    \end{align*}
    Thus, \citet[Corollary 1]{GaObRe2015} guarantees the existence of a barrier function $r^{\nu,\eta}$ such that, if $\,W_0\sim\frac{\nu}{\nu[\mathbb{R}]}$, i.e.\ distributed like the normalized measure of $\nu$, then $\,W_\varsigma\sim\frac{\eta}{\eta[\mathbb{R}]}\,$ where $\,\varsigma=\inf\{t\geq0\,|\,t\geq r^{\nu,\eta}(W_t)\}$.
    
    3) By construction, $\,\eta[(\frac{3a+b}{4},\frac{a+3b}{4})]=0$, because that is true for all $\xi_c$. Hence, as all two-sided hitting times of Brownian motion are unbounded and hit both barriers, if $\,\nu[(\frac{3a+b}{4},\frac{a+3b}{4})]>0$, then $\,r^{\nu,\eta}(x)=\infty\,$ for all $\,x\in(\frac{3a+b}{4},\frac{a+3b}{4})$.
    
    4) As $\,\nu[\mathbb{R}\!\setminus\!(a,b)]=0$, Lemma \ref{lem:mu[tail]=0} ensures that $\,r^{\nu,\eta}(x)=0\,$ for all $\,x\notin(a,b)$.
    
    Now, we define the claimed measure $\mu$ implicitly through its barrier function $r$. Fix a number $x\in\mathbb{R}$ and choose a sequence $\,x_i\downarrow x\,$ for $\,i\uparrow\infty$, let $\,\nu_i=\Phi|_{(x_{i+1},x_i)}\,$ be the standard normal distribution on $(x_{i+1},x_i)$ for $\,i\in\mathbb{N}$, and define
    \[
        r(y)=
        \begin{cases}
            1+r^{\nu_i,\eta}(y)&\mbox{if $\,y\in(x_{i+1},x_i)$}
            \\
            \hfill1&\mbox{if $\,y\in\mathbb{R}\setminus\bigcup_{i=1}^\infty(x_{i+1},x_i)$}
        \end{cases}.
    \]    
    To emphasize the dependence of $\eta$ on $\nu$, we write $\,\eta=\eta_\nu$. By construction, see 2) and 4), the probability measure $\mu$ corresponding to $r$, i.e.\ the image measure of $W_\tau$ with $\,\tau=\inf\{t\geq0\,|\,t\geq r(W_t)\}$, is given by
    \[
        \mu=
        \begin{cases}
            \eta_{\nu_i}&\mbox{on $\,(x_{i+1},x_i)$}
            \\
            \hfill\Phi&\mbox{on $\,\mathbb{R}\setminus\bigcup_{i=1}^\infty(x_{i+1},x_i)$}
        \end{cases}
    .\]
    Note, $r$ is indeed a barrier function of a Root solution, because $\,t\mapsto W_{t\wedge\tau}$ is a uniformly integrable martingale. To see this, observe that large values of $W_{t\wedge\tau}$ imply that $\,W_{t\wedge\tau}=W_{t\wedge1}$, and $\,t\mapsto W_{t\wedge1}\,$ is a uniformly integrable martingale, in particular,
    \[\lim_{K\rightarrow\infty}\sup_{t\geq0}\,\mathbb{E}\big[|W_{t\wedge\tau}|\mathbbm{1}_{W_{t\wedge\tau}>K}\big]
    =\lim_{K\rightarrow\infty}\sup_{t\geq0}\,\mathbb{E}\big[|W_{t\wedge1}|\mathbbm{1}_{W_{t\wedge1}>K}\big]
    =0.\]
    Moreover, $r$ is regular, because $\,u_{\mu_0}(z)\neq u_\mu(z)\,$ for all $z$. To see this, observe that Jensen's inequality is strict for $\,y\mapsto|z-y|\,$ with $\,z\in\mathbb{R}$, because $\,\mu=\Phi\,$ on $\,\mathbb{R}\setminus[x,x_1]$. Also, as $\,t\mapsto W_{t\wedge \tau}\,$ is a uniformly integrable martingale, it holds true that $\,0=\mathbb{E}[W_0]=\mathbb{E}[W_\tau]$. Altogether,
    \[\int_\mathbb{R}\!|z-y|\,\mathrm{d}\mu(y)
    >\Big|\int_\mathbb{R}\!z-y\,\mathrm{d}\mu(y)\Big|
    =\big|z-\mathbb{E}[W_\tau]\big|
    =|z|
    =\int_\mathbb{R}\!|z-y|\,\mathrm{d}\mu_0(y).\]
    Furthermore, 1) implies that $\mu$ has a density and 3) implies that $\,r(\frac{x_i+x_{i+1}}{2})=\infty\,$ for all $i$.
    
    Overall, $\mu$ is an absolutely continuous probability measure with $\,\mu[(\frac{3x_{i+1}+x_i}{4},\frac{x_{i+1}+3x_i}{4})]=0\,$ and the barrier function $r$ corresponding to Root's solution fulfills $\,r(x_i)=1\,$ and $\,r(\frac{x_i+x_{i+1}}{2})=\infty\,$ for all $i$.
\end{proof}
\end{ex}

\bibliographystyle{abbrvnat} 
\bibliography{Stochastics}
\end{document}